\documentclass[12pt]{article} 

\usepackage[utf8]{inputenc} 

\usepackage{geometry} 
\usepackage{graphicx} 

\usepackage{booktabs} 
\usepackage{array} 
\usepackage{verbatim} 
\usepackage{subfig} 

\usepackage{fancyhdr} 
\usepackage{ifthen}

\newcommand{ \randbin }[1] {X_{#1}}  
\newcommand{ \countuniq }[1] {
  \sigma_{#1}
} 
\newcommand{ \randcount }[2] {
  C^{#1}_{#2}
} 
\newtheorem{theorem}{Theorem}

\begin{document}

\title{The Number of Distinct Subsequences of a Random Binary String}
\author{Michael J. Collins\\ Christopher Newport University\\Department of Physics, Computer Science, and Engineering\\\texttt{michael.collins@cnu.edu}}

\maketitle
\begin{abstract}
  We determine the average number of distinct subsequences in a random binary string, and derive an estimate for the average number of distinct subsequences of a particular length.
\end{abstract}

\section{Introduction}
Let $\randbin{n}$ be a uniformly distributed random binary string of length $n$. For any string $s$, let $\countuniq{m}(s)$ be the number of \emph{distinct} subsequences of $s$ of length $m$, and let $\countuniq{}(s)$ be the total number of distinct subsequences. We consider the random variables $\randcount{n}{m} = \countuniq{m}(\randbin{n})$ and $\randcount{n}{} = \countuniq{}(\randbin{n})$. Note that we are counting sub\emph{sequences}, not substrings (the latter would have to be contiguous portions of $s$). Questions about such random varaibles arise in the study of error-correcting codes for deletion channels \cite{MR2525669} -- communication channels in which each bit has some probability of being dropped, with the receiver getting no indication of where these deletions occurred. When a codeword $W$ is sent on such a channel, the receiver sees a subsequence of $W$.

We denote the expectation of a random variable $A$ by $\hat A$. Given a string $s = (s_0s_1\cdots s_{n-1})$, we write $s[i:j]$ for the substring $s = (s_is_{i+1}\cdots s_{j-1})$. A \emph{run} is a maximal constant substring. A sequence of length $m$ is called an $m$-sequence (similarly an $m$-subseqence, $m$-string et cetera).

\subsection{Subsequences of length $m$}
If a string $s$ begins with a run of length $k$, then we can assume with no loss of generality -- since we are only counting \emph{distinct} subsequences -- that  every nonempty subsequence starts at either $s_0$ or $s_k$. Now $\countuniq{m}(s)$ equals the number of distinct $m$-subsequences that start with zero, plus the number that start with 1, i.e.
\begin{equation}
\countuniq{m}(s) = \countuniq{m-1}(s[1:n]) + \countuniq{m-1}(s[k+1:n]).
\end{equation} 
where the second term is zero when $k=n$.
Of course $\randbin{n}$ has probability $2^{-k}$ of starting with a run of length $k$ if $k<n$, and probability $2^{1-n}$ of starting with a run of length $n$ (i.e. being a constant sequence); thus (by linearity of expectation) $\hat\randcount{n}{m}$ satisfies the recurrence
\begin{equation}
  \label{eq:countuniqM}
\hat \randcount{n}{m} = \hat \randcount{n-1}{m-1} + \sum_{k \geq 1}^{n-1} 2^{-k} \hat\randcount{n-k-1}{m-1}  
\end{equation}
Since $\hat\randcount{n}{m}=0$ when $m>n$ we can (for $m>1$) write this as
\begin{equation}
  \label{eq:countuniqM2}
\hat \randcount{n}{m}  = \sum_{k \geq 0}^{\infty} 2^{-k} \hat\randcount{n-k-1}{m-1}
\end{equation}
from which we obtain a recurrence relation similar to the recurrence for the binomial coefficents:
\begin{theorem}
\begin{equation}
\hat \randcount{n}{m} = \hat \randcount{n-1}{m-1} + \frac{1}{2}\hat\randcount{n-1}{m}
\end{equation}
with initial conditions $\hat\randcount{n}{n}=1$ and $\hat\randcount{n}{0}=1$.
\end{theorem}

Note in particular that $\hat\randcount{n}{n-1} = \frac{n+1}{2}$. This must be the expected number of runs in $X_n$, since an $(n-1)$-subsequence is determined entirely by the run from which one bit is deleted. Indeed the first bit of $\randbin{n}$ starts a new run with probability 1, while each subsequent bit starts a new run with probability $\frac{1}{2}$, again giving (by linearity of expectation) $\hat\randcount{n}{n-1}=\frac{n+1}{2}$.

More generally,  $\hat\randcount{n}{n-m}$ is, for fixed $m$, a polynomial of degree $m$. 
Let $\hat\randcount{n}{n-m} = p_m(n) = \sum_{i=0}^m \alpha_{m,i} n^i$. Now
\[
p_m(n) = p_m(n-1) + \frac{1}{2}p_{m-1}(n-1)
\]
thus equating coefficients on $n^{m-1}$ gives
\[
\alpha_{m,m-1} = -m\alpha_{m,m} + \alpha_{m,m-1} + \frac{1}{2}\alpha_{m-1,m-1}
\]
so
\[
\alpha_{m,m} = \frac{1}{2m}\alpha_{m-1,m-1} \ .
\]
Thus (since $\alpha_{0,0} = 1$) we obtain $\alpha_m = 1/(2^mm!)$ and we have the approximation
\begin{theorem}
\begin{equation}
  \label{eq:binomialApprox}
  \hat\randcount{n}{n-m} = 2^{-m} {n \choose m} + O(n^{m-1})
\end{equation}
\end{theorem}


\subsection{Total Number of Subsequences}
By the same reasoning as in (\ref{eq:countuniqM}) we have
\begin{equation}
\countuniq{}(s) = 1+ \countuniq{}(s[1:n]) + \countuniq{}(s[k+1:n]).
\end{equation}
when $s$ begins with a $k$-run (the initial $1$ counts the empty subsequence), and thus 
\begin{equation}
  \label{eq:countuniqAll}
  \hat\randcount{n}{} = 1 + \hat\randcount{n-1}{} + \sum_{k \geq 1}2^{-k}\hat\randcount{n-k-1}{}
  = 1 + \sum_{k \geq 0}2^{-k}\hat\randcount{n-k-1}{} \ .
\end{equation}
Noting that
\begin{equation}
  \sum_{k \geq 1}2^{-k}\hat\randcount{n-k-1}{} = \frac{\hat\randcount{n-1}{}-1}{2}
\end{equation}
we have
\begin{equation}
  \hat\randcount{n}{} = \frac{1}{2} + \frac{3}{2}\hat\randcount{n-1}{}
\end{equation}
and $\hat\randcount{0}{}=1$, thus
\begin{theorem}
\begin{equation}
  \label{eq:countuniqAllSltn}
  \hat\randcount{n}{} = 2\left( \frac{3}{2} \right)^n - 1 \ .
\end{equation}
\end{theorem}
This improves an earlier result \cite{FlaxmanHS04} that $\hat\randcount{n}{} = O((3/2)^n)$.

\bibliographystyle{plain}
\bibliography{distinctSubseq}

\end{document}